\documentclass{amsart}
\usepackage{amsfonts,latexsym,amsmath,amscd,amssymb}
\usepackage{bbm}
\usepackage{a4wide}

\theoremstyle{plain}
\newtheorem{theorem}{Theorem}
\numberwithin{theorem}{section}

\newtheorem{corollary}{Corollary}
\numberwithin{corollary}{section}

\numberwithin{definition}{section}

\newtheorem{lemma}{Lemma}
\numberwithin{lemma}{section}

\newtheorem{proposition}{Proposition}
\numberwithin{proposition}{section}

\numberwithin{remark}{section}

 \numberwithin{equation}{section}

\newcommand {\be}{\begin{equation}}
\newcommand {\ee}{\end{equation}}

\newcommand{\h}{\begin{eqnarray*}}
 \newcommand{\e}{\end{eqnarray*}}
\newcommand{\CC}{\mathbb{C}}
\begin{document}
\title[Witten Genus and String Complete Intersections]{Witten Genus and String Complete Intersections}

\author{Qingtao Chen}
\address{Q. Chen, Department of Mathematics, University of California,
Berkeley, CA, 94720-3840} \email{chenqtao@math.berkeley.edu}
\date{Jan 11, 2007}

\author{Fei Han}
\address{F. Han, \ Department of Mathematics, University of California,
Berkeley, CA, 94720-3840} \email{feihan@math.berkeley.edu}


\maketitle

\begin{abstract}We prove that the Witten genus of certain nonsingular string complete intersections in products of
complex projective spaces vanishes. Our result generalizes a known
result of Landweber and Stong (cf. [HBJ]).
\end{abstract}

\section {Introduction} Let $M$ be a $4k$ dimensional closed
oriented smooth manifold. Let $E$ be a complex vector bundle over
$M$. For any complex number $t$, set
$\Lambda_t(E)=\CC|M+tE+t^2\Lambda^2(E)+\cdots$ and
$S_t(E)=\CC|M+tE+t^2S^2(E)+\cdots$, where for any integer $j\geq
1, \Lambda^j(E)$ (resp. $S^j(E)$) is the $j$-th exterior (resp.
symmetric) power of $E$ (cf. [At]). Set
$\widetilde{E}=E-\CC^{\mathrm{rk}(E)}$.

Let $q=e^{\pi i \tau}$ with $\tau\in \mathbb{H}$, the upper half
plane. Define after Witten ([W]) \be \Theta_q(E)=\bigotimes_{n\geq
1}S_{q^{2n}}(E).\ee The Witten genus ([W]) $\varphi_W(M)$ for $M$
is defined as \be
\varphi_W(M)=\langle\widehat{A}(M)\mathrm{ch}(\Theta_q(\widetilde{TM\otimes
\CC})), [M]\rangle,\ee where $\widehat{A}(M)$ is the
$\widehat{A}$-characteristic class of $M$ and $[M]$ is the
fundamental class of $M$. Let $\{\pm 2\pi \sqrt{-1} x_j, 1\leq j
\leq 2k\} $ be the formal Chern roots of $TM\otimes \CC$. Then the
Witten genus can be written by using Chern roots as (cf.
[L2][L3])\be
\varphi_W(M)=\langle\prod_{j=1}^{2k}x_j\frac{\theta'(0,
\tau)}{\theta(x_j, \tau)}, [M]\rangle,\ee where $\theta(v, \tau)$
is the Jacobi theta function (see (2.1) below). The Witten genus
was first introduced by E. Witten [W] in studying quantum filed
theory and can be viewed as the loop space analogue of the
$\widehat{A}$ genus.

According to the Atiyah-Singer index theorem, when $M$ is spin,
$\varphi_W(M)\in \mathbb{Z}[[q]]$ (cf. [HBJ]). Moreover, when the
spin manifold $M$ is string, i.e. $\frac{p_1(TM)}{2}=0$
($\frac{p_1(TM)}{2}$ is a dimension 4 characteristic class, twice
of which is the first integral Pontryagin class $p_1(TM)$), or
even weaker, when the first rational Pontryagin class $p_1(M)=0$,
$\varphi_W(M)$ is a modular form of weight $2k$ with integer
Fourier expansion (cf. [HBJ]).

Let $V_{(d_{pq})}$ be a nonsingular $4k$ dimensional complete
intersection in the product of complex projective spaces ${\CC
P}^{n_1}\times {\CC P}^{n_2}\times \cdots \times {\CC P}^{n_s}$,
which is dual to $\prod_{p=1}^t(\sum_{q=1}^sd_{pq}x_q)\in
H^{2t}({\CC P}^{n_1}\times {\CC P}^{n_2}\times \cdots \times {\CC
P}^{n_s}, \mathbb{Z})$, where $x_q\in H^2({\CC P}^{n_q},
\mathbb{Z}), 1\leq q \leq s,$ is the generator of $H^*({\CC
P}^{n_q}, \mathbb{Z})$ and $d_{pq}, 1\leq p\leq t, 1\leq q \leq
s,$ are integers. Let $P_q: {\CC P}^{n_1}\times {\CC
P}^{n_2}\times \cdots \times {\CC P}^{n_s}\rightarrow {\CC
P}^{n_q}, 1\leq q \leq s,$ be the $q$-th projection. Then
$V_{(d_{pq})}$ is the intersection of the zero loci of smooth
global sections of line bundles
$\otimes_{q=1}^{s}P_q^*(\mathcal{O}(d_{pq})), 1 \leq p \leq t$,
where $\mathcal{O}(d_{pq})=\mathcal{O}(1)^{d_{pq}}$ is the
$d_{pq}$-th power of the canonical line bundle $\mathcal{O}(1)$
over ${\CC P}^{n_q}$.

We should point out that we have abused the terminology, complete
intersection, in algebraic geometry. Actually, since we don't
require that the integers $d_{pq}$'s be nonnegative,
$V_{(d_{pq})}$ might not be an algebraic variety. However, by
transversality, $V_{(d_{pq})}$ can always be chosen to be smooth.
Putting some relevant conditions (see Proposition 3.1 below) on
the data $n_q, 1 \leq q \leq s$ and $d_{pq}, 1\leq p\leq t, 1\leq
q \leq s$, the complete intersection $V_{(d_{pq})}$ can be made
string. This systematically provides us a lot of interesting
examples of string manifolds. This paper is devoted to study the
Witten genus of string manifolds generated in this way. See also
[GM] and [GO] for a study of elliptic genera of complete
intersections and the Landau-Ginzburg/CalabiYau correspondence.

Let $$D=\left[\begin{array}{cccc}
                                      d_{11}&d_{12}& \cdots &d_{1s}\\
                                      d_{21}&d_{22}&\cdots &d_{2s}\\
                                      \cdots &\cdots & \cdots   &\cdots \\
                                      d_{t1}&d_{t2}&\cdots&d_{ts}
                                     \end{array}\right].$$
Let $m_q$ be the number of nonzero elements in the $q$-th column
of $D$.

The main result of this paper is
\begin{theorem} If $m_q+2 \leq n_q, 1\leq q \leq s$ and $V_{(d_{pq})}$ is string, then the Witten genus
$\varphi_W(V_{(d_{pq})})$ vanishes.
\end{theorem}

Our result generalizes a known result that the Witten genus of
string nonsingular complete intersections of hypersurfaces with
degrees $d_1, \cdots, d_t$ in a single complex projective space
vanishes, which is due to Landweber and Stong (cf. [HBJ, Sect.
6.3]). In [HBJ], the proof of the result of Landweber and Stong is
roughly described by applying the properties of the sigma
function. Also in this special case, our result is broader than
Landweber-Stong's result, since we don't require $d_1, \cdots,
d_t$ be all positive.

Explicitly expanding $\Theta_q(\widetilde{TM\otimes \CC})$, we get
\be\widehat{A}(M)\mathrm{ch}(\Theta_q(\widetilde{TM\otimes
\CC}))=\widehat{A}(M)+\widehat{A}(M)(\mathrm{ch}(TM\otimes
\CC)-4k)q^2+\cdots.\ee Therefore it's not hard to obtain the
following corollary from Theorem 1.1.

\begin{corollary} If $m_q+2 \leq n_q, 1\leq q \leq s$ and $V_{(d_{pq})}$ is string, then $$\langle
\widehat{A}(V_{(d_{pq})}), [V_{(d_{pq})}]\rangle=0, \ \langle
\widehat{A}(V_{(d_{pq})})\mathrm{ch}(TV_{(d_{pq})}\otimes \CC),
[V_{(d_{pq})}]\rangle=0.$$

\end{corollary}

Let $M$ be a 12 dimensional oriented closed smooth manifold. The
signature of $M$ can be expressed by the $\widehat{A}$-genus and
the twisted $\widehat{A}$-genus as the following [AGW][L1], \be
L(M)=8\widehat{A}(M)\mathrm{ch}(T_\CC M)-32\widehat{A}(M). \ee
Combining Corollary 1.1 and (1.5), we obtain that

\begin{corollary} If $m_q+2 \leq n_q, 1\leq q \leq s$, $V_{(d_{pq})}$ is 12-dimensional and
string, then the signature of $V_{(d_{pq})}$ vanishes.
\end{corollary}

Let $M$ be a 16 dimensional oriented closed smooth manifold. One
has the following formula [CH], \be L(M)\mathrm{ch}(T_\CC
M)=-2048\{\widehat{A}(M)\mathrm{ch}(T_\CC
M)-48\widehat{A}(M)\}.\ee Combining Corollary 1.1 and (1.6), for
twisted signature $\mathrm{Sig}(M, T)\triangleq \langle
L(TM)\mathrm{ch}(T_\CC M), [M]\rangle$, we obtain that
\begin{corollary} If $m_q+2 \leq n_q, 1\leq q \leq s$, $V_{(d_{pq})}$ is 16-dimensional and
string, then the twisted signature  $\mathrm{Sig}(V_{(d_{pq})},
T)$ vanishes.
\end{corollary}

In the following, we will review some necessary preliminaries in
Section 2 and prove Theorem 1.1 in Section 3.

\section{Some preliminaries}
In this section, we review some tools and knowledge that we are
going to apply in the proof of Theorem 1.1.

First, let's review the relevant concepts and results on residues
in complex geometry. See chapter 5 of [GH] for details.

Let $U$ be the ball $\{z\in {\CC}^s:\|z\|<\varepsilon\}$ and $f_1,
\cdots, f_s\in \mathcal{O}(\overline{U})$ functions holomorphic in
a neighborhood of the closure $\overline{U}$ of $U$. We assume
that the $f_i(z)$ have the origin as isolated common zero. Set
$$D_i=(f_i)=\mathrm{divisor\ of} f_i,$$
$$D=D_1+\cdots+D_s.$$
Let $$\omega=\frac{g(z)dz_1\wedge \cdots \wedge dz_s}{f_1(z)\cdots
f_s(z)}$$ be a meromorphic $s$-form with polar divisor $D$. The
{\bf residue} of $\omega$ at the origin is defined as the
following
$$ \mathrm{Res}_{\{0\}}\omega=\left(\frac{1}{2\pi \sqrt{-1}} \right)^s\int_\Gamma
\omega,$$ where $\Gamma$ is the real $s$-cycle defined by
$$\Gamma=\{z:|f_i(z)|=\varepsilon, 1 \leq i \leq s\}$$
and oriented by
$$d(\mathrm{arg}f_1)\wedge\cdots \wedge d(\mathrm{arg}
f_s)\geq0.$$

Let $M$ be a compact complex manifold of dimension $s$. Suppose
that $D_1, \cdots, D_s$ are effective divisors, the intersection
of which is a finite set of points. Let $D=D_1+\cdots+D_s.$ Let
$\omega$ be a meromorphic $s$-form on $M$ with polar divisor $D$.
For each point $P\in D_1\cap\cdots \cap D_s$, we may restrict
$\omega$ to a neighborhood of $U_P$ of $P$ and define the residue
$\mathrm{Res}_P\omega$ as above. One has (cf. [GH], Chaper 5)
\begin{lemma}(\bf Residue Theorem)
$$\Sigma_{P\in  D_1\cap\cdots \cap D_s}\mathrm{Res}_P\omega=0.$$
\end{lemma}

We also need some knowledge on the Jacobi theta functions.
Although we are going to use only one of them, for the sake of
completeness, we list the definitions and transformation laws of
all of them.

The four Jacobi theta functions are defined as follows (cf. [Ch]):
\be\theta(v,\tau)=2q^{1/4}\sin(\pi v)
\prod_{j=1}^\infty\left[(1-q^{2j})(1-e^{2\pi
\sqrt{-1}v}q^{2j})(1-e^{-2\pi \sqrt{-1}v}q^{2j})\right]\ ,\ee \be
\theta_1(v,\tau)=2q^{1/4}\cos(\pi v)
 \prod_{j=1}^\infty\left[(1-q^{2j})(1+e^{2\pi \sqrt{-1}v}q^{2j})
 (1+e^{-2\pi \sqrt{-1}v}q^{2j})\right]\ ,\ee
\be \theta_2(v,\tau)=\prod_{j=1}^\infty\left[(1-q^{2j})
 (1-e^{2\pi \sqrt{-1}v}q^{2j-1})(1-e^{-2\pi \sqrt{-1}v}q^{2j-1})\right]\
 ,\ee
\be \theta_3(v,\tau)=\prod_{j=1}^\infty\left[(1-q^{2j}) (1+e^{2\pi
\sqrt{-1}v}q^{2j-1})(1+e^{-2\pi \sqrt{-1}v}q^{2j-1})\right]\ ,\ee
where $q=e^{\pi i \tau}$, $\tau \in \mathbb{H}$, the upper half
plane and $v\in \CC$. They are all holomorphic functions for
$(v,\tau)\in \mathbb{C \times H}$. Let
$\theta^{'}(0,\tau)=\frac{\partial}{\partial
v}\theta(v,\tau)|_{v=0}$. They satisfy the following relations
(cf. [Ch]): \be \theta(v+1, \tau)=-\theta(v, \tau),\
\theta(v+\tau, \tau)=-{1\over q}e^{-2\pi i v}\theta(v, \tau),\ee
\be \theta_1(v+1, \tau)=-\theta_1(v, \tau),\  \theta_1(v+\tau,
\tau)={1\over q}e^{-2\pi i v}\theta_1(v, \tau),\ee \be
\theta_2(v+1, \tau)=\theta_2(v, \tau),\  \theta_2(v+\tau,
\tau)=-{1\over q}e^{-2\pi i v}\theta_2(v, \tau),\ee \be
\theta_3(v+1, \tau)=\theta_3(v, \tau),\  \theta_3(v+\tau,
\tau)={1\over q}e^{-2\pi i v}\theta_3(v, \tau).\ee Therefore it's
not hard to deduce in the following how the theta functions vary
along the lattice $\Gamma=\{m+n\tau| m,n\in \mathbb{Z}\}$. We have
\be \theta(v+m, \tau)=(-1)^m\theta(v, \tau)\ee and \be
\begin{split} &\theta(v+n\tau, \tau)\\
=&-{1\over q}e^{-2\pi i(v+(n-1)\tau)}\theta(v+(n-1)\tau, \tau)\\
=&-{1\over q}e^{-2\pi i(v+(n-1)\tau)}\left(-{1\over
q}\right)e^{-2\pi
i(v+(n-2)\tau)}\theta(v+(n-2)\tau, \tau)\\
=&(-1)^n\frac{1}{q^n}e^{-2\pi i[(v+(n-1)\tau)+(v+(n-2)\tau)+\cdots
+v]}\theta(v,\tau)\\
=&(-1)^n\frac{1}{q^n}e^{-2\pi i nv-\pi i
n(n-1)\tau}\theta(v,\tau)\\
=&(-1)^ne^{-2\pi inv-\pi i n^2\tau}\theta(v,\tau).
\end{split} \ee Similarly, we have
\be \theta_1(v+m, \tau)=(-1)^m\theta_1(v, \tau), \
\theta_1(v+n\tau, \tau)=e^{-2\pi inv-\pi i
n^2\tau}\theta_1(v,\tau);\ee  \be \theta_2(v+m, \tau)=\theta_2(v,
\tau), \ \theta_2(v+n\tau, \tau)=(-1)^ne^{-2\pi inv-\pi i
n^2\tau}\theta_2(v,\tau);\ee \be \theta_3(v+m, \tau)=\theta_3(v,
\tau), \ \theta_3(v+n\tau, \tau)=e^{-2\pi inv-\pi i
n^2\tau}\theta_3(v,\tau).\ee

\section{Proof of Theorem 1.1} Let $i:V_{(d_{pq})}\rightarrow {\CC P}^{n_1}\times {\CC
P}^{n_2}\times \cdots \times {\CC P}^{n_s}$ be the inclusion. It's
not hard to see that \be i^* T_\mathbb{R}({\CC P}^{n_1}\times {\CC
P}^{n_2}\times \cdots \times {\CC P}^{n_s})\cong
TV_{(d_{pq})}\oplus
i^*\left(\oplus_{p=1}^{t}(\otimes_{q=1}^{s}P_q^*
\mathcal{O}(d_{pq}))\right),\ee where we forget the complex
structure of the line bundles $\otimes_{q=1}^{s}P_q^*
\mathcal{O}(d_{pq}), 1\leq p \leq t.$ Therefore for the total
Stiefel-Whitney class, we have\be i^*w( T_\mathbb{R}({\CC
P}^{n_1}\times {\CC P}^{n_2}\times \cdots \times {\CC P}^{n_s}))=
w(TV_{(d_{pq})})\prod_{p=1}^{t}i^*w(\otimes_{q=1}^{s}P_q^*
\mathcal{O}(d_{pq})),\ee or more precisely \be
i^*\left(\prod_{q=1}^{s}(1+x_q)^{n_q+1}\right)\equiv
w(TV_{(d_{pq})})\prod_{p=1}^{t}i^*\left(1+\sum_{q=1}^{s}d_{pq}x_q\right)\
\ \  \ \ \ \mathrm{mod}\, 2. \ee By (3.3), we can easily see that
\be w_1(TV_{(d_{pq})})=0, \ee \be w_2(TV_{(d_{pq})})\equiv
\sum_{q=1}^{s}\left(n_q+1-\sum_{p=1}^{t}d_{pq}\right)i^*x_q\ \ \ \
\mathrm{mod}\,2.\ee

As for the total rational Pontryagin class, we have \be i^*p(
T_\mathbb{R}({\CC P}^{n_1}\times {\CC P}^{n_2}\times \cdots \times
{\CC P}^{n_s}))=
p(TV_{(d_{pq})})\prod_{p=1}^{t}i^*p(\otimes_{q=1}^{s}P_q^*
\mathcal{O}(d_{pq})),\ee or \be
p(V_{(d_{pq})})=\prod_{q=1}^s(1+(i^*x_q)^2)^{n_q+1}\prod_{p=1}^t\left(1+\left(\sum_{q=1}^sd_{pq}i^*x_q\right)^2\right)^{-1}.\ee
Hence we have \be \begin{split}
p_1(V_{(d_{pq})})&=\sum_{q=1}^{s}(n_q+1)(i^*x_q)^2-\sum_{p=1}^{t}\left(\sum_{q=1}^sd_{pq}i^*x_q\right)^2\\
&=\sum_{q=1}^s(n_q+1-\sum_{p=1}^{t}d_{pq}^2)(i^*x_q)^2-\sum_{1\leq
u, v \leq s, u\neq
v}\left(\sum_{p=1}^{t}d_{pu}d_{pv}i^*x_ui^*x_v\right).\end{split}\ee

Let $i_!: H^*(V_{(d_{pq})}, \mathbb{Q})\rightarrow H^{*+2t}({\CC
P}^{n_1}\times {\CC P}^{n_2}\times \cdots \times {\CC P}^{n_s},
\mathbb{Q})$ be the push forward. Thus if $p_1(V_{(d_{pq})})=0$,
then
$$i_!p_1(V_{(d_{pq})})=i_!i^*\left(\sum_{q=1}^s(n_q+1-\sum_{p=1}^{t}d_{pq}^2)x_q^2-\sum_{1\leq
u, v \leq s, u\neq
v}\left(\sum_{p=1}^{t}d_{pu}d_{pv}x_ux_v\right)\right)=0,$$ i.e.
$$\left(\prod_{p=1}^t\left(\sum_{q=1}^sd_{pq}x_q\right)\right)\left(\sum_{q=1}^s(n_q+1-\sum_{p=1}^{t}d_{pq}^2)x_q^2-\sum_{1\leq
u, v \leq s, u\neq
v}\left(\sum_{p=1}^{t}d_{pu}d_{pv}x_ux_v\right)\right)=0
$$ in $ H^{2t+4}({\CC P}^{n_1}\times {\CC P}^{n_2}\times \cdots
\times {\CC P}^{n_s}, \mathbb{Q}). $ If $m_q+2 \leq n_q, 1\leq q
\leq s$, then the left hand side of the above equality should not
only be a zero element in the cohomology ring but also be a zero
polynomial. Note that the polynomial ring is an integral domain.
Therefore at least one of it's factor should be zero. But
$\prod_{p=1}^t\left(\sum_{q=1}^sd_{pq}x_q\right)$ is nonzero. This
means
$$\sum_{q=1}^s(n_q+1-\sum_{p=1}^{t}d_{pq}^2)x_q^2-\sum_{1\leq u, v
\leq s, u\neq v}(\sum_{p=1}^{t}d_{pu}d_{pv}x_ux_v)=0$$ and
consequently the following identities hold, \be
n_q+1-\sum_{p=1}^{t}d_{pq}^2=0, \ 1\leq q\leq s;\ee \be
\sum_{p=1}^{t}d_{pu}d_{pv}=0,\ 1\leq u, v \leq s, u\neq v.\ee

Note that \be n_q+1-\sum_{p=1}^{t}d_{pq}^2\equiv
n_q+1-\sum_{p=1}^{t}d_{pq}\ \ \ \\ \mathrm{mod}\,2,\ \ \ \  1\leq
q\leq s. \ee Hence (3.9) implies that $w_2(TV_{(d_{pq})})=0$.

In a summary, we have the following proposition.
\begin{proposition} When $m_q+2 \leq n_q, 1\leq q
\leq s$, $p_1(V_{(d_{pq})})=0$ implies $V_{(d_{pq})}$ is spin.
Therefore when $m_q+2 \leq n_q, 1\leq q \leq s$, $V_{(d_{pq})}$ is
string if and only if one of the following holds \newline (1)
$p_1(V_{(d_{pq})})=0$; \newline (2) the following identities hold,
\be n_q+1-\sum_{p=1}^{t}d_{pq}^2=0, \ 1\leq q\leq s, \ee \be
\sum_{p=1}^{t}d_{pu}d_{pv}=0,\ 1\leq u, v \leq s, u\neq v;\ee
\newline (3) in the matrix
$$D=\left[\begin{array}{cccc}
                                      d_{11}&d_{12}& \cdots &d_{1s}\\
                                      d_{21}&d_{22}&\cdots &d_{2s}\\
                                      \cdots &\cdots & \cdots   &\cdots \\
                                      d_{t1}&d_{t2}&\cdots&d_{ts}
                                     \end{array}\right]$$
$\parallel\mathrm{col_qD}\parallel^2=n_q+1, 1 \leq q \leq s$ and
any two columns are orthogonal to each other; or equivalently,
$$D^tD=\mathrm{diag}(n_1+1, \cdots, n_s+1). $$

\end{proposition}

With the above preparations, we are able to prove Theorem 1.1.
$$ $$
\noindent{\it Proof of Theorem 1.1:} Let $[V_{(d_{pq})}]$ be the
fundamental class of $V_{(d_{pq})}$ in $H_{4k}(V_{(d_{pq})},
\mathbb{Z})$. Then according to (1.3) and the multiplicative
property of the Witten genus, up to a constant scalar,
\be \begin{split}&\varphi_W(V_{(d_{pq})})\\
=&\left(\left(\prod_{q=1}^{s}\left[\frac{i^*x_q}{\frac{\theta(i^*x_q,\tau)}{\theta'(0,
\tau)}}\right]^{n_q+1}\right)\left(\prod_{p=1}^{t}\left[\frac{\sum_{q=1}^sd_{pq}i^*x_q}{\frac{\theta(\sum_{q=1}^sd_{pq}i^*x_q,\tau)}{\theta'(0,
\tau)}}\right]^{-1}\right)\right)[V_{(d_{pq})}]\\
=&\left(\left(\prod_{q=1}^{s}\left[\frac{x_q}{\frac{\theta(x_q,\tau)}{\theta'(0,
\tau)}}\right]^{n_q+1}\right)\left(\prod_{p=1}^{t}\left[\frac{1}{\frac{\theta(\sum_{q=1}^sd_{pq}x_q,\tau)}{\theta'(0,
\tau)}}\right]^{-1}\right)\right)[{\CC P}^{n_1}\times {\CC
P}^{n_2}\times \cdots \times {\CC P}^{n_s}]\\
=&\mathrm{coefficient\ of}\ x_1^{n_1}\cdots x_s^{n_s} \mathrm{in}\
\left(\frac{(\prod_{q=1}^{s}{x_q}^{n_q+1})\left(\prod_{p=1}^{t}\frac{\theta(\sum_{q=1}^sd_{pq}x_q,\tau)}{\theta'(0,
\tau)}\right)}{\prod_{q=1}^{s}\left[\frac{\theta(x_q,\tau)}{\theta'(0,
\tau)}\right]^{n_q+1}}\right)\\
=&\mathrm{Res}_{0}\left(\frac{\left(\prod_{p=1}^{t}\frac{\theta(\sum_{q=1}^sd_{pq}x_q,\tau)}{\theta'(0,
\tau)}\right)dx_1\wedge\cdots \wedge
dx_s}{\prod_{q=1}^{s}\left[\frac{\theta(x_q,\tau)}{\theta'(0,
\tau)}\right]^{n_q+1}}\right).\\
\end{split}\ee
Note that in (3.14), we have used Pioncare duality to deduce the
second equality.

Set $$g(x_1, \cdots,
x_s)=\prod_{p=1}^{t}\frac{\theta(\sum_{q=1}^sd_{pq}x_q,\tau)}{\theta'(0,
\tau)}, \ \ \ \ \ \
f_q(x_q)=\left[\frac{\theta(x_q,\tau)}{\theta'(0,
\tau)}\right]^{n_q+1}, 1\leq q \leq s,$$ and $$\omega=\frac{g(x_1,
\cdots, x_s) dx_1\wedge\cdots \wedge dx_s}{f_1(x_1)\cdots
f_s(x_s)}.$$ Then up to a constant scalar, \be
\varphi_W(V_{(d_{pq})})=\mathrm{Res}_{(0, 0, \cdots, 0)}\omega.\ee

By (2.9), \be \begin{split} &g(x_1+1, x_2, \cdots,
x_s)\\
=&\prod_{p=1}^{t}\frac{\theta(\sum_{q=1}^sd_{pq}x_q+d_{p1},\tau)}{\theta'(0,
\tau)}\\
=&(-1)^{d_{11}+\cdots+d_{t1}}\prod_{p=1}^{t}\frac{\theta(\sum_{q=1}^sd_{pq}x_q,\tau)}{\theta'(0,
\tau)}\end{split}\ee and
$$f_1(x_1+1)=\left[\frac{\theta(x_1+1,\tau)}{\theta'(0,
\tau)}\right]^{n_1+1}=(-1)^{n_1+1}\left[\frac{\theta(x_1,\tau)}{\theta'(0,
\tau)}\right]^{n_1+1}.$$ Thus $\frac{g(x_1+1, \cdots, x_s)
}{f_1(x_1+1)\cdots
f_s(x_s)}=(-1)^{(d_{11}+\cdots+d_{t1})-(n_1+1)}\frac{g(x_1,
\cdots, x_s) }{f_1(x_1)\cdots f_s(x_s)}.$ Note that by (3.12),
$$(d_{11}+\cdots+d_{t1})-(n_1+1)\equiv (d_{11}^2+\cdots+d_{t1}^2)-(n_1+1)=0\
\mathrm{mod}\,2.$$ Thus one obtains that $\frac{g(x_1+1, \cdots,
x_s) }{f_1(x_1+1)\cdots f_s(x_s)}=\frac{g(x_1, \cdots, x_s)
}{f_1(x_1)\cdots f_s(x_s)}.$ Similarly, we have \be \frac{g(x_1,
\cdots,x_q+1, \cdots, x_s) }{f_1(x_1)\cdots f_q(x_q+1)\cdots
f_s(x_s)}=\frac{g(x_1, \cdots, x_s) }{f_1(x_1)\cdots f_s(x_s)},
1\leq q \leq s.\ee

On the other hand, by (2.9) and (2.10), \be
\begin{split}&g(x_1+\tau, x_2, \cdots,
x_s)\\=&\prod_{p=1}^{t}\frac{\theta(\sum_{q=1}^sd_{pq}x_q+d_{p1}\tau,\tau)}{\theta'(0,
\tau)}\\
=&\prod_{p=1}^{t}(-1)^{d_{p1}}e^{-2\pi i
d_{p1}(\sum_{q=1}^sd_{pq}x_q)-\pi i
d_{p1}^2\tau}\frac{\theta(\sum_{q=1}^sd_{pq}x_q,\tau)}{\theta'(0,
\tau)}\\
=&(-1)^{d_{11}+\cdots+d_{t1}}e^{-2\pi
i\sum_{p=1}^{t}d_{p1}(\sum_{q=1}^sd_{pq}x_q)-\pi i\tau
(\sum_{p=1}^{t}d_{p1}^2)}\frac{\theta(\sum_{q=1}^sd_{pq}x_q,\tau)}{\theta'(0,
\tau)}\\
=&(-1)^{d_{11}+\cdots+d_{t1}}e^{-2\pi
i\sum_{p=1}^{t}d_{p1}(\sum_{q=1}^sd_{pq}x_q)-\pi i\tau
(\sum_{p=1}^{t}d_{p1}^2)}g(x_1, x_2, \cdots, x_s)\\
\end{split}\ee
and \be \begin{split}
&f_1(x_1+\tau)\\
=&\left[\frac{\theta(x_1+\tau,\tau)}{\theta'(0,
\tau)}\right]^{n_1+1}\\
=&\left[-e^{-2\pi i x_1-\pi
i\tau}\frac{\theta(x_1,\tau)}{\theta'(0,
\tau)}\right]^{n_1+1}\\
=&(-1)^{n_1+1}e^{-2\pi i (n_1+1)x_1-\pi
i\tau(n_1+1)}\left[\frac{\theta(x_1,\tau)}{\theta'(0,
\tau)}\right]^{n_1+1}\\
=&(-1)^{n_1+1}e^{-2\pi i (n_1+1)x_1-\pi
i\tau(n_1+1)}f_1(x_1).\end{split}\ee Therefore \be
\begin{split}&\frac{g(x_1+\tau, \cdots, x_s) }{f_1(x_1+\tau)\cdots
f_s(x_s)}\\
=& (-1)^{d_{11}+\cdots+d_{t1}-n_1-1}e^{-2\pi
i\sum_{p=1}^{t}d_{p1}(\sum_{q=1}^sd_{pq}x_q)-\pi i\tau
(\sum_{p=1}^{t}d_{p1}^2)+2\pi i (n_1+1)x_1+\pi
i\tau(n_1+1)}\\
&\cdot \frac{g(x_1, \cdots, x_s) }{f_1(x_1)\cdots
f_s(x_s)}.\end{split}\ee However \be
d_{11}+\cdots+d_{t1}-n_1-1\equiv d_{11}^2+\cdots+d_{t1}^2-n_1-1 \
\mathrm{mod} 2 \ee and \be \begin{split}&-2\pi
i\sum_{p=1}^{t}d_{p1}\left(\sum_{q=1}^sd_{pq}x_q\right)-\pi i\tau
(\sum_{p=1}^{t}d_{p1}^2)+2\pi i (n_1+1)x_1+\pi i\tau(n_1+1)\\
=&\pi i \tau \left[(n_1+1)-\sum_{p=1}^{t}d_{p1}^2\right]+2\pi
i\left[(n_1+1)-\sum_{p=1}^{t}d_{p1}^2\right]x_1-2\pi
i\sum_{q=2}^{s}\left(\sum_{p=1}^{t}d_{p1}d_{pq}x_q\right).\end{split}\ee
Therefore by (3.12) and (3.13),
$$ d_{11}+\cdots+d_{t1}-n_1-1\equiv 0\ \mathrm{mod}\,2,$$
and $$-2\pi
i\sum_{p=1}^{t}d_{p1}\left(\sum_{q=1}^sd_{pq}x_q\right)-\pi i\tau
(\sum_{p=1}^{t}d_{p1}^2)+2\pi i (n_1+1)x_1+\pi i\tau(n_1+1)=0.$$
Consequently, by (3.20), we obtain that \be\frac{g(x_1+\tau,
\cdots, x_s) }{f_1(x_1+\tau)\cdots f_s(x_s)}=\frac{g(x_1, \cdots,
x_s) }{f_1(x_1)\cdots f_s(x_s)}.\ee

Similarly, one also obtains that \be \frac{g(x_1, \cdots,x_q+\tau,
\cdots, x_s) }{f_1(x_1)\cdots f_q(x_q+\tau)\cdots
f_s(x_s)}=\frac{g(x_1, \cdots, x_s) }{f_1(x_1)\cdots f_s(x_s)}, \
1\leq q \leq s. \ee Therefore from (3.17) and (3.24), we see that
$\omega$ can be viewed as a meromorphic $s$-form defined on the
product of $s$-tori, $\left(\CC/\Gamma\right)^s$, which is a
compact complex manifold.

$\theta(v, \tau)$ has the lattice points $m+n\tau, m,n\in
\mathbb{Z}$ as it's simple zero points [Ch]. We therefore see that
$\omega$ has pole divisors $\{0\}\times
\left(\CC/\Gamma\right)^{s-1}, \left(\CC/\Gamma\right)\times
\{0\}\times \left(\CC/\Gamma\right)^{s-2}, \cdots,
\left(\CC/\Gamma\right)^{s-1}\times \{0\}$. So $(0, 0, \cdots, 0)$
is the unique intersection point of these polar divisors.

Therefore by the residue theorem on compact complex manifolds, we
directly deduces that $\mathrm{Res}_{(0, 0, \cdots, 0)}\omega=0$.
By (3.15), we obtain that
$\varphi_W(V_{(d_{pq})})=\mathrm{Res}_{(0, 0, \cdots,
0)}\omega=0$. Q.E.D

\section{Acknowledgement} F. Han is grateful to Professor Peter Teichner for a lot of discussions and helps.
He also thanks Professor Kefeng Liu for inspiring suggestions. Q.
Chen is grateful to Professor Nicolai Reshetikhin for his interest
and support. Thanks also go to Professor Friedrich Hirzebruch,
Professor Stephan Stolz and Professor Weiping Zhang for their
interests and many discussions with us. We would like to thank
Professor Michael Joachim and Professor Serge Ochanine for
inspiring communications with us. The paper was finished when the
second author was visiting the Max-Planck-Institut
f$\ddot{\mathrm{u}}$r Mathematik at Bonn.

\end{document}